\documentclass[12pt]{article}

\usepackage{amsmath, epsfig, cite}
\usepackage{amssymb}
\usepackage{amsfonts}
\usepackage{latexsym}
\usepackage{amsthm}

\makeatletter
\renewcommand{\@seccntformat}[1]{{\csname the#1\endcsname}.\hspace{.5em}}
\makeatother

\newtheorem{thm}{Theorem}[section]

\newtheorem{cor}[thm]{Corollary}
\newtheorem{conj}[thm]{Conjecture}
\newtheorem{lem}[thm]{Lemma}

\newcommand{\pf}{\noindent{\it Proof.} }

\def\Q{\mathbb Q}

\renewcommand{\qed}{\hfill$\Box$\medskip}

\numberwithin{equation}{section}

\begin{document}

\renewcommand{\thefootnote}{*}

\begin{center}
{\Large\bf Some $q$-analogues of supercongruences of Rodriguez-Villegas\\[5pt] }
\end{center}

\vskip 2mm \centerline{Victor J. W. Guo$^1$  and Jiang Zeng$^{2}$}
\begin{center}
{\footnotesize $^1$Department of Mathematics, Shanghai Key Laboratory of PMMP,
East China Normal University,\\ 500 Dongchuan Rd., Shanghai 200241,
 People's Republic of China\\
{\tt jwguo@math.ecnu.edu.cn,\quad http://math.ecnu.edu.cn/\textasciitilde{jwguo}}\\[10pt]
$^2$Universit\'e de Lyon; Universit\'e Lyon 1; Institut Camille
Jordan, UMR 5208 du CNRS;\\ 43, boulevard du 11 novembre 1918,
F-69622 Villeurbanne Cedex, France\\
{\tt zeng@math.univ-lyon1.fr,\quad
http://math.univ-lyon1.fr/\textasciitilde{zeng}} }
\end{center}

\vskip 0.7cm \noindent{\small{\bf Abstract.}}
 We study different $q$-analogues and generalizations of
the ex-conjectures of Rodriguez-Villegas. For example, for any odd
prime $p$, we show that the  known congruence
$$
\sum_{k=0}^{p-1}\frac{{2k\choose k}^2}{16^k} \equiv
\left(\frac{-1}{p}\right)\pmod{p^2},
$$
where $(\frac{\cdot}{p})$ is the Legendre symbol, has the following two nice $q$-analogues:
\begin{align*}
\sum_{k=0}^{p-1}\frac{(q;q^2)_k^2 }{(q^2;q^2)_k^2}q^{(1+\varepsilon)k}
&\equiv \left(\frac{-1}{p}\right) q^{\frac{(p^2-1)\varepsilon}{4}}\pmod{(1+q+\cdots+q^{p-1})^2}, 
\end{align*}
where $(a;q)_n=(1-a)(1-aq)\cdots(1-aq^{n-1})$ and $\varepsilon=\pm1$. Several related conjectures are also proposed.



\section{Introduction}

Rodriguez-Villegas \cite{RV} discovered numerically some remarkable
supercongruences between a truncated hypergeometric function associated to a Calabi-Yau
manifold at a prime $p$ and the number of its $\mathbb{F}_p$-points. In particular,
Rodriguez-Villegas recorded four such supercongruences associated to elliptic curves.
Following a strategy
 developed by Ahlgren and Ono \cite{AO}, by
using  the Gross-Koblitz formula to write the Gaussian hypergeometric
series in terms of the $p$-adic $\Gamma$-function,
Mortenson \cite{Mortenson1,Mortenson4} first proved  the following  four conjectured
supercongruences  of Rodriguez-Villegas\cite[(36)]{RV}.

 \begin{thm} [Rodriguez-Villegas-Mortenson]\label{thm:Mortenson}
Let $p\geqslant 5$ be a prime. Then
\begin{align}
\sum_{k=0}^{p-1}\frac{{2k\choose k}^2}{16^k}
&\equiv \left(\frac{-1}{p}\right)\pmod{p^2}, \label{eq:RV1} \\
\sum_{k=0}^{p-1}\frac{{3k\choose 2k}{2k\choose k}}{27^k}
&\equiv \left(\frac{-3}{p}\right)\pmod{p^2}, \label{eq:RV2} \\
\sum_{k=0}^{p-1}\frac{{4k\choose 2k}{2k\choose k}}{64^k}
&\equiv \left(\frac{-2}{p}\right)\pmod{p^2}, \label{eq:RV3} \\
\sum_{k=0}^{p-1}\frac{{6k\choose 3k}{3k\choose k}}{432^k}
&\equiv \left(\frac{-1}{p}\right)\pmod{p^2}, \label{eq:RV4}
\end{align}
where $(\frac{\cdot}{p})$ denotes the Legendre symbol modulo $p$.
\end{thm}
Elementary proof of Theorem~\ref{thm:Mortenson}  has been given by Z.-H. Sun~\cite{SunZH2}. See also
\cite{CLZ,SunZH1,Tauraso,Tauraso1} for several simple proofs of \eqref{eq:RV1}.
A generalization of \eqref{eq:RV1} to the modulus $p^3$ case was obtained by Z.-W. Sun \cite{Sun4}. Note that van Hamme \cite{Hamme} and
McCarthy and Osburn \cite{MO} have studied some related interesting supercongruences.

Recall that the {\it $q$-shifted factorials} are defined by $(a;q)_0=1$ and
$$(a;q)_n=(1-a)(1-aq)\cdots(1-aq^{n-1})\ \text{for}\ n=1,2,\ldots,$$
and the {\it $q$-integer} is defined as $[p]= 1+q +\cdots+ q^{p-1}$. The starting point of this paper is the observation of the following striking $q$-analogue of
Theorem~\ref{thm:Mortenson}.
\begin{conj}\label{conj:2346}
Let $p\geqslant 5$ be a prime and let $(\frac{\cdot}{p})$
be the Legendre symbol modulo $p$. Then
\begin{align*}
\sum_{k=0}^{p-1}\frac{(q;q^2)_k^2}{(q^2;q^2)_k^2}
&\equiv \left(\frac{-1}{p}\right)q^{\frac{1-p^2}{4}}\pmod{[p]^2}, \\
\sum_{k=0}^{p-1}\frac{(q;q^3)_k (q^2;q^3)_k}{(q^3;q^3)_k^2}
&\equiv \left(\frac{-3}{p}\right)q^{\frac{1-p^2}{3}}\pmod{[p]^2}, \\
\sum_{k=0}^{p-1}\frac{(q;q^4)_k (q^3;q^4)_k}{(q^4;q^4)_k^2}
&\equiv \left(\frac{-2}{p}\right)q^{\frac{3(1-p^2)}{8}}\pmod{[p]^2}, \\
\sum_{k=0}^{p-1}\frac{(q;q^6)_k (q^5;q^6)_k}{(q^6;q^6)_k^2}
&\equiv \left(\frac{-1}{p}\right)q^{\frac{5(1-p^2)}{12}}\pmod{[p]^2}.
\end{align*}
\end{conj}

Congruences modulo $[p]$ or $[p]^2$ (even $[p]^3$) have been studied by different authors
(see \cite{Andrews99,Pan,SP,Straub,Tauraso2}).
Throughout the paper we will tacitly use the fact that when $p$ is a prime the $q$-integer
$[p]$ is an irreducible polynomial in $\Q[q]$.
Therefore $\Q[q]/[p]$ is a field.
Furthermore, rational functions $t(q)/s(q)$ are
well defined modulo $[p]$ or $[p]^r$ ($r\geqslant 1$) provided that $s(q)$ is relatively prime to $[p]$.
For any two polynomials
$A(x,q)=\sum_{k=0}^{n} a_k (q)x^k$ and $B(x,q)=\sum_{k=0}^{n} b_k(q)x^k$ in $x$
with coefficients being rational functions  $t(q)/s(q)$  such that $s(q)$ is relatively prime to $[p]$,
we  use the convention that
$$
A(x,q)\equiv B(x,q)\pmod{[p]^r}\Longleftrightarrow a_k(q)\equiv b_k(q) \pmod{[p]^r}\quad \textrm{for}
\quad k=0, 1,\ldots,n.
$$

There are  several generalizations and variations of \eqref{eq:RV1}--\eqref{eq:RV4} in the literature, but
no $q$-analogues seem to be investigated hitherto.
Indeed, Tauraso \cite{Tauraso}
proved the following  generalization of \eqref{eq:RV1}.
\begin{thm}[Tauraso \cite{Tauraso}]Let $p$ be an odd prime. Then
\begin{align}
\sum_{k=0}^{p-1}{2k\choose k}^2\frac{x^k}{16^k}
&\equiv
\sum_{k=0}^{\frac{p-1}{2}}{\frac{p-1}{2}\choose k}^2 (-x)^k (1-x)^{\frac{p-1}{2}-k} \pmod{p^2}. \label{eq:Tauraso}
\end{align}
\end{thm}

Recently, Z.-H. Sun \cite{SunZH2} introduced the generalized Legendre polynomials
\begin{align}
P_n(a,x)=\sum_{k=0}^n{a\choose k}{-1-a\choose k}\frac{(1-x)^k}{2^k}
=\sum_{k=0}^n{a\choose k}{a+k\choose k}\frac{(x-1)^k}{2^k}, \label{eq:def-sun}
\end{align}
and proved many supercongruences related to $P_{p-1}(a,x)$. In particular, he obtained
the following result.
\begin{thm} [Z.-H. Sun \cite{SunZH2}]
Let $p$ be an odd prime and let $a$ be a $p$-adic integer. Then
\begin{align}
P_{p-1}(a,x)\equiv (-1)^{\langle a \rangle_p} P_{p-1}(a,-x)\pmod{p^2}, \label{eq:sun-lengendre-1}
\end{align}
and so
\begin{align}
\sum_{k=0}^{p-1}{a\choose k}{-1-a\choose k}(x^k-(-1)^{\langle a\rangle_p} (1-x)^k)
\equiv 0 \pmod{p^2},  \label{eq:sun-lengendre-2}
\end{align}
where $\langle a \rangle_p$ denotes the least nonnegative residue of $a$ modulo $p$.
\end{thm}
It is easy to see that the Rodriguez-Villegas-Mortenson congruences \eqref{eq:RV1}--\eqref{eq:RV4}
immediately follows from the congruence \eqref{eq:sun-lengendre-2} by taking $x=1$
and $a=-\frac{1}{2},-\frac{1}{3},-\frac{1}{4},-\frac{1}{6}$.
The aim of this paper is to give $q$-analogues of \eqref{eq:Tauraso}--\eqref{eq:sun-lengendre-2}.
It turns out that a complete $q$-analogue of \eqref{eq:Tauraso} is easily given. However, for a general
$p$-adic integer $a$, we can only give $q$-analogues of \eqref{eq:sun-lengendre-1} and
\eqref{eq:sun-lengendre-2} in the modulus $p$ case. On the other hand, for $a=-\frac{1}{2}$, we shall give
complete $q$-analogues of them.
Thus, the first congruence in Conjecture \ref{conj:2346} is proved, while the other
three congruences are still open. Some further related unsolved problems will also be presented in this paper.

\section{Results, I: supercongruences modulo $[p]^2$}
Recall that the {\it $q$-binomial coefficients} ${n\brack k}$ are defined by
$$
{n\brack k}={n\brack k}_q
=\begin{cases}
\displaystyle\frac{(q^{n-k+1};q)_k}{(q;q)_k}, &\text{if $0\leqslant k\leqslant n$},\\[10pt]
0,&\text{otherwise.}
\end{cases}
$$
We first give a $q$-analogue of Tauraso's congruence \eqref{eq:Tauraso}.
\begin{thm}\label{thm:q-Tauraso}
Let $p$ be an odd prime. Then
\begin{align}
\sum_{k=0}^{p-1}\frac{(q;q^2)_k^2}{(q^2;q^2)_k^2} x^k
\equiv
\sum_{k=0}^{\frac{p-1}{2}}{\frac{p-1}{2}\brack k}_{q^2}^2 q^{-pk+k^2} (-x)^k (x;q^2)_{\frac{p-1}{2}-k}\pmod{[p]^2}.
\label{eq:q-Tauraso}
\end{align}
\end{thm}
Since
$$
\lim_{q\to 1}\frac{(q;q^2)_k^2}{(q^2;q^2)_k^2}
=\left(\prod_{j=1}^{k}\frac{2j-1}{2j}\right)^2={2k\choose k}^216^{-k},
$$
letting $q\to 1$ in \eqref{eq:q-Tauraso}, we obtain
 \eqref{eq:Tauraso}.
 Moreover, setting $x=1$ in \eqref{eq:q-Tauraso} yields  the following $q$-analogue of \eqref{eq:RV1}.
\begin{cor}\label{cor:2.2}Let $p$ be an odd prime. Then
\begin{align}
\sum_{k=0}^{p-1}\frac{(q;q^2)_k^2}{(q^2;q^2)_k^2}
&\equiv \left(\frac{-1}{p}\right)q^{\frac{1-p^2}{4}}\pmod{[p]^2}. \label{eq:q-RV1}
\end{align}
\end{cor}

\noindent{\it Remark.} Corollary \ref{cor:2.2} confirms the first congruence in Conjecture \ref{conj:2346}.
\medskip

Our second result is another generalization of \eqref{eq:q-RV1}.
\begin{thm}\label{thm:vq-Tauraso}
Let $p$ be an odd prime. Then
\begin{align}
\sum_{k=0}^{p-1}\frac{(q;q^2)_k^2}{(q^2;q^2)_k^2}x^k
&\equiv \left(\frac{-1}{p}\right) q^{\frac{1-p^2}{4}}
\sum_{k=0}^{p-1}\frac{(q;q^2)_k^2}{(q^2;q^2)_k^2} q^{2k}(x;q^2)_k \pmod{[p]^2}. \label{eq:q2-RV1}
\end{align}
\end{thm}
It is clear that, when $x=1$, the congruence \eqref{eq:q2-RV1} reduces to \eqref{eq:q-RV1}. On the other
hand, setting $x=0$ in \eqref{eq:q2-RV1}, we obtain the following dual form of \eqref{eq:q-RV1}.
\begin{cor}Let $p$ be an odd prime. Then
\begin{align*}
\sum_{k=0}^{p-1}\frac{(q;q^2)_k^2 }{(q^2;q^2)_k^2}q^{2k}
&\equiv \left(\frac{-1}{p}\right)q^{\frac{p^2-1}{4}}\pmod{[p]^2}. 
\end{align*}
\end{cor}

Our third result is a $q$-analogue of the $a=-\frac{1}{2}$ case of \eqref{eq:sun-lengendre-1}.
\begin{thm}\label{thm:p-121}
Let $p$ be an odd prime. Then
\begin{align}
\sum_{k=0}^{p-1}\frac{(q;q^2)_k^2 (x;q^2)_k q^{2k}}{(q^2;q^2)_k^2 (-q^2;q^2)_k}
\equiv \left(\frac{-1}{p}\right)
\sum_{k=0}^{p-1}\frac{(q;q^2)_k^2 (-x;q^2)_k q^{2k}}{(q^2;q^2)_k^2 (-q^2;q^2)_k}
\pmod{[p]^2}. \label{eq:thm-main2-00}
\end{align}
\end{thm}
Letting $x=-1$ in \eqref{eq:thm-main2-00}, and noticing that $\frac{(-1;q^2)_k}{(-q^2;q^2)_k}=\frac{2}{1+q^{2k}}$,
we are led to another $q$-analogue of \eqref{eq:RV1}.
\begin{cor}Let $p$ be an odd prime. Then
\begin{align}
\sum_{k=0}^{p-1}\frac{2(q;q^2)_k^2 q^{2k}}{(q^2;q^2)_k^2 (1+q^{2k})}
&\equiv \left(\frac{-1}{p}\right)\pmod{[p]^2}.
\end{align}
\end{cor}
If $\left(\frac{-1}{p}\right)=-1$ and $x=0$, then we immediately deduce that
both sides of \eqref{eq:thm-main2-00} are congruent to $0$ modulo $[p]^2$, which may be restated as follows.
\begin{cor}Let $p$ be a prime of the form $4k+3$. Then
\begin{align}
\sum_{k=0}^{p-1}\frac{(q;q^2)_k^2 q^{2k}}{(q^2;q^2)_k^2 (-q^2;q^2)_k}
\equiv 0  \pmod{[p]^2}. \label{eq:qBeukers}
\end{align}
\end{cor}
Note that, when $q=1$, the congruence \eqref{eq:qBeukers} can be written as
\begin{align*}
\sum_{k=0}^{p-1}\frac{{2k\choose k}^2}{32^k}\equiv 0  \pmod{p^2}\quad\text{for}\quad p\equiv 3\pmod 4,
\end{align*}
which was conjecture by Z.-W. Sun \cite{Sun0} and proved by Tauraso \cite{Tauraso} and
Z.-H. Sun \cite{SunZH1,SunZH2}.

\section{Results, II: congruences modulo $[p]$}
In this section, we first  give $q$-analogues of \eqref{eq:RV2}--\eqref{eq:RV4}.
Actually we shall prove the following  more general results.
\begin{thm}\label{thm:2mp2}
Let $p$ be an odd prime and $m$, $r$ two positive integers with $p\nmid m$.  Then
\begin{align}
&\hskip -2mm \sum_{k=0}^{p-1}\frac{(q^r;q^m)_k (q^{m-r};q^m)_k x^k}{(q^m;q^m)_k^2} \notag\\
&\equiv (-1)^{\langle -\frac{r}{m}\rangle_p}
q^{\frac{-m\langle -\frac{r}{m}\rangle_p \left(\langle -\frac{r}{m}\rangle_p+1 \right)}{2}}
\sum_{k=0}^{p-1}\frac{(q^r;q^m)_k (q^{m-r};q^m)_k (x;q^m)_k q^{mk}}{(q^m;q^m)_k^2} \pmod{[p]}. \label{eq:2mp2}
\end{align}
In particular, if $p\equiv \pm 1\pmod m$, then
\begin{align}
&\hskip -2mm \sum_{k=0}^{p-1}\frac{(q^r;q^m)_k (q^{m-r};q^m)_k x^k}{(q^m;q^m)_k^2} \notag\\
&\equiv (-1)^{\langle -\frac{r}{m}\rangle_p}q^{\frac{r(m-r)(1-p^2)}{2m}}
\sum_{k=0}^{p-1}\frac{(q^r;q^m)_k (q^{m-r};q^m)_k (x;q^m)_k q^{mk}}{(q^m;q^m)_k^2} \pmod{[p]}. \label{eq:2mp2-2}
\end{align}
\end{thm}
Note that, for $p\geqslant 5$,  we have
\begin{align*}
(-1)^{\langle -\frac{1}{3}\rangle_p}
&=\begin{cases}
(-1)^{\frac{p-1}{3}}=1,&\text{if $p\equiv 1\pmod 3$} \\
(-1)^{\frac{2p-1}{3}}=-1,&\text{if $p\equiv 2\pmod 3$}
\end{cases}
=\left(\frac{-3}{p}\right), \\
(-1)^{\langle -\frac{1}{4}\rangle_p}
&=\begin{cases}
(-1)^{\frac{p-1}{4}},&\text{if $p\equiv 1\pmod 4$} \\
(-1)^{\frac{3p-1}{4}},&\text{if $p\equiv 3\pmod 4$}
\end{cases}
=\left(\frac{-2}{p}\right), \\
(-1)^{\langle -\frac{1}{6}\rangle_p}
&=\begin{cases}
(-1)^{\frac{p-1}{6}},&\text{if $p\equiv 1\pmod 6$} \\
(-1)^{\frac{5p-1}{6}},&\text{if $p\equiv 5\pmod 6$}
\end{cases}
=(-1)^{\frac{p-1}{2}}=\left(\frac{-1}{p}\right).
\end{align*}
Letting $x=1$, $r=1$, $m=3,4,6$ in \eqref{eq:2mp2-2} with $p\geqslant 5$,
we obtain the following result.
\begin{cor}\label{cor:one}
Let $p\geqslant 5$ be a prime. Then
\begin{align*}
\sum_{k=0}^{p-1}\frac{(q;q^3)_k (q^2;q^3)_k}{(q^3;q^3)_k^2}
&\equiv \left(\frac{-3}{p}\right)q^{\frac{1-p^2}{3}}\pmod{[p]}, \\
\sum_{k=0}^{p-1}\frac{(q;q^4)_k (q^3;q^4)_k}{(q^4;q^4)_k^2}
&\equiv \left(\frac{-2}{p}\right)q^{\frac{3(1-p^2)}{8}}\pmod{[p]}, \\
\sum_{k=0}^{p-1}\frac{(q;q^6)_k (q^5;q^6)_k}{(q^6;q^6)_k^2}
&\equiv \left(\frac{-1}{p}\right)q^{\frac{5(1-p^2)}{12}}\pmod{[p]}.
\end{align*}
\end{cor}
\noindent{\it Remark.} The congruences in Corollary \ref{cor:one} confirm
the remaining three congruences in Conjecture \ref{conj:2346} modulo $[p]$.
\medskip

In the same vein, letting $x=0$, $r=1$, and $m=3,4,6$ in \eqref{eq:2mp2}, we obtain
\begin{cor}\label{cor:two}
Let $p\geqslant 5$ be a prime. Then
\begin{align*}
\sum_{k=0}^{p-1}\frac{(q;q^3)_k (q^2;q^3)_k}{(q^3;q^3)_k^2}q^{3k}
&\equiv \left(\frac{-3}{p}\right)q^{\frac{p^2-1}{3}}\pmod{[p]}, \\
\sum_{k=0}^{p-1}\frac{(q;q^4)_k (q^3;q^4)_k}{(q^4;q^4)_k^2}q^{4k}
&\equiv \left(\frac{-2}{p}\right)q^{\frac{3(p^2-1)}{8}}\pmod{[p]}, \\
\sum_{k=0}^{p-1}\frac{(q;q^6)_k (q^5;q^6)_k}{(q^6;q^6)_k^2}q^{6k}
&\equiv \left(\frac{-1}{p}\right)q^{\frac{5(p^2-1)}{12}}\pmod{[p]}.
\end{align*}
\end{cor}

It seems that we have the following stronger result (see Conjecture \ref{conj:one} for a further generalization).
\begin{conj}The congruences in Corollary {\rm\ref{cor:two}} hold modulo $[p]^2$.
\end{conj}

\begin{thm}\label{thm:2mp2-new}
Let $p$ be an odd prime and $m$, $r$ two positive integers with $p\nmid m$.
Then
\begin{align*}
&\hskip -2mm \sum_{k=0}^{p-1}\frac{(q^r;q^m)_k (q^{m-r};q^m)_k x^k}{(q^m;q^m)_k^2} \\
&\equiv
\sum_{k=0}^{\langle -\frac{r}{m}\rangle_p}{\langle -\frac{r}{m}\rangle_p\brack k}_{q^m}^2
q^{\frac{mk(k-1)}{2}-mk\langle -\frac{r}{m}\rangle_p} (-x)^k (x;q^m)_{\langle -\frac{r}{m}\rangle_p-k} \pmod{[p]}.
\end{align*}
\end{thm}

Next, we give  $q$-analogues of
\eqref{eq:sun-lengendre-1}--\eqref{eq:sun-lengendre-2} in the modulus $p$ case. For this end,
we introduce the following  polynomial
$$
P_{n,m,r}(q,x)
=\sum_{k=0}^{n}\frac{(q^r;q^m)_k (q^{m-r};q^m)_k (x;q^m)_k q^{mk}}{(q^m;q^m)_k^2 (-q^m;q^m)_k}.
$$
Note  that $P_{n,m,r}(1,x)$ is  the generalized
Legendre polynomial $P_{n}(a,x)$ in \eqref{eq:def-sun} with $a=-\frac{r}{m}$.
\begin{thm}\label{thm:main2}
Let $p$ be an odd prime and $m$, $r$ two positive integers with $p\nmid m$.  Then
\begin{align}
P_{p-1,m,r}(q,x)
\equiv (-1)^{\langle -\frac{r}{m}\rangle_p} P_{p-1,m,r}(q,-x) \pmod{[p]}. \label{eq:thm-main2}
\end{align}
\end{thm}
Letting $x=0$ in \eqref{eq:thm-main2}, we obtain
\begin{cor}Let $p$ be an odd prime and $m,r$ two integers with $p\nmid m$ and
${\langle -r/m\rangle_p}\equiv 1\pmod 2$. Then
\begin{align}
\sum_{k=0}^{p-1}\frac{(q^r;q^m)_k (q^{m-r};q^m)_k  q^{mk}}{(q^m;q^m)_k^2 (-q^m;q^m)_k}
\equiv 0  \pmod{[p]}.  \label{eq:amp=1}
\end{align}
\end{cor}
Taking $(m,r)=(3,1),(4,1),(6,1)$ in \eqref{eq:amp=1}, we get the following congruences.
\begin{cor}\label{cor:three}
Let $p$ be an odd prime. Then
\begin{align}
\sum_{k=0}^{p-1}\frac{(q;q^3)_k (q^2;q^3)_k q^{3k}}{(q^3;q^3)_k^2 (-q^3;q^3)_k}
&\equiv 0  \pmod{[p]},\quad\text{for}\quad p\equiv 2\pmod 3,\\[5pt]
\sum_{k=0}^{p-1}\frac{(q;q^4)_k (q^3;q^4)_k q^{4k}}{(q^4;q^4)_k^2 (-q^4;q^4)_k}
&\equiv 0 \pmod{[p]},\quad\text{for}\quad p\equiv 5,7\pmod 8,\\[5pt]
\sum_{k=0}^{p-1}\frac{(q;q^6)_k (q^5;q^6)_k q^{6k}}{(q^6;q^6)_k^2 (-q^6;q^6)_k}
&\equiv 0 \pmod{[p]},\quad\text{for}\quad p\equiv 3\pmod 4.
\end{align}
\end{cor}
Letting $x=-1$ in \eqref{eq:thm-main2}, we obtain
\begin{align}
P_{p-1,m,r}(q,-1)
\equiv (-1)^{\langle -\frac{r}{m}\rangle_p}  \pmod{[p]}.
\end{align}
Taking $(m,r)=(3,1),(4,1),(6,1)$, we get the following result.
\begin{cor}\label{cor:four}
Let $p\geqslant 5$ be a prime. Then
\begin{align}
\sum_{k=0}^{p-1}\frac{2(q;q^3)_k (q^2;q^3)_k q^{3k}}{(q^3;q^3)_k^2 (1+q^{3k})}
&\equiv \left(\frac{-3}{p}\right)\pmod{[p]}, \\
\sum_{k=0}^{p-1}\frac{2(q;q^4)_k (q^3;q^4)_k q^{4k}}{(q^4;q^4)_k^2 (1+q^{4k})}
&\equiv \left(\frac{-2}{p}\right)\pmod{[p]}, \\
\sum_{k=0}^{p-1}\frac{2(q;q^6)_k (q^5;q^6)_k q^{6k}}{(q^6;q^6)_k^2 (1+q^{6k})}
&\equiv \left(\frac{-1}{p}\right)\pmod{[p]}.
\end{align}
\end{cor}

\section{Proofs of Theorems \ref{thm:q-Tauraso} and \ref{thm:vq-Tauraso}}
Recall that the little $q$-Legendre polynomials are defined by
\begin{align}
P_n(x|q)=\sum_{k=0}^{n}{n\brack k}{n+k\brack k} q^{\frac{k(k+1)}{2}-nk} (-x)^k.  \label{eq:pnxq-1}
\end{align}
They can also be written as (see \cite{VAssche})
\begin{align}
P_n(x|q)=(-1)^n q^{-\frac{n(n+1)}{2}}\sum_{k=0}^{n}{n\brack k}{n+k\brack k}(-1)^k q^{\frac{k(k+1)}{2}-nk} (xq;q)_k.
\label{eq:pnxq-2}
\end{align}
We now give a new expansion for the little $q$-Legendre polynomials.
\begin{lem}Let $n$ be a nonnegative integer. Then
\begin{align}
P_n(x|q)=\sum_{k=0}^{n}{n\brack k}^2 q^{\frac{k(k+1)}{2}-nk} (-x)^k (xq;q)_{n-k}. \label{eq:pnxq-3}
\end{align}
\end{lem}
\pf  By the $q$-binomial theorem (see \cite[p.~36, Theorem 3.3]{Andrews}):
\begin{align}
(x;q)_N=\sum_{k=0}^N {N\brack k}(-x)^k q^{\frac{k(k-1)}{2}}, \label{eq:qbinomial}
\end{align}
one sees that, for $0\leqslant m\leqslant n$,
the coefficient of $x^m$ in the right-hand side of \eqref{eq:pnxq-3} is given by
\begin{align*}
&\hskip -2mm
\sum_{k=0}^{m}{n\brack k}^2 q^{\frac{k(k+1)}{2}-nk} (-1)^k {n-k\brack m-k}(-1)^{m-k}q^{\frac{(m-k)(m-k+1)}{2}} \\[5pt]
&=(-1)^m {n\brack m}\sum_{k=0}^{m}{m\brack k}{n\brack k} q^{(m-k)(n-k)-mn+\frac{m(m+1)}{2}} \\[5pt]
&=(-1)^m {n\brack m}{n+m\brack m}q^{-mn+\frac{m(m+1)}{2}},
\end{align*}
where the last step follows from the $q$-Chu-Vandermonde formula (see \cite[p.~37, Theorem 3.4]{Andrews}).
This completes the proof.
\qed

We also need the following result.
\begin{lem}\label{lem:one}
Let $p$ be an odd prime and $0\leqslant k\leqslant p-1$. Then
\begin{align}
\frac{(q;q^2)_k^2}{(q^2;q^2)_k^2}
\equiv (-1)^k {\frac{p-1}{2}\brack k}_{q^2} {\frac{p-1}{2}+k\brack k}_{q^2} q^{k^2-kp} \pmod{[p]^2}.
\label{eq:lem1}
\end{align}
\end{lem}
\pf Observing that
\begin{align*}
(1-q^{2j-1})^2+(1-q^{p-2j+1})(1-q^{p+2j-1})q^{2j-1-p}
=(1-q^p)^2 q^{2j-1-p},
\end{align*}
we have
\begin{align*}
(1-q^{2j-1})^2
\equiv
-(1-q^{p-2j+1})(1-q^{p+2j-1})q^{2j-1-p} \pmod{[p]^2}.
\end{align*}
It follows that
\begin{align*}
\frac{(q;q^2)_k^2}{(q^2;q^2)_k^2}
=\prod_{j=1}^k \frac{(1-q^{2j-1})^2}{(1-q^{2j})^2}
&\equiv (-1)^k\prod_{j=1}^k \frac{(1-q^{p-2j+1})(1-q^{p+2j-1})q^{2j-1-p}}{(1-q^{2j})^2} \\
&= (-1)^k {\frac{p-1}{2}\brack k}_{q^2} {\frac{p-1}{2}+k\brack k}_{q^2} q^{k^2-kp} \pmod{[p]^2},
\end{align*}
as desired.  \qed

\noindent{\it Proof of Theorem {\rm\ref{thm:q-Tauraso}.}}
Letting $n=\frac{p-1}{2}$, replacing $q$ and $x$ by $q^2$ and $xq^{-2}$ respectively
in \eqref{eq:pnxq-1} and \eqref{eq:pnxq-3},
and then applying \eqref{eq:lem1}, we obtain \eqref{eq:q-Tauraso}.
\qed

\noindent{\it Proof of Theorem {\rm\ref{thm:vq-Tauraso}.}}
The proof is similar to that of Theorem \ref{thm:q-Tauraso} by just comparing \eqref{eq:pnxq-1} and \eqref{eq:pnxq-2}.
\qed

\section{Proof of Theorem \ref{thm:p-121}}
We fist establish two lemmas.
\begin{lem}\label{lem:andrews}
Let $n$ be a positive integer and $0\leqslant j\leqslant n$. Then
\begin{align}
&\hskip -2mm
\sum_{k=j}^{n}(-1)^k{n+k\brack k}{n-j\brack k-j}\frac{q^{\frac{k(k+1)}{2}-nk}}{(-q;q)_k} \notag\\
&=\begin{cases}
\displaystyle (-1)^{\frac{n-j}{2}}{n\brack \frac{n}{2}}_{q^2}
\frac{(q^{n+1};q^2)_{\frac{j}{2}}\, q^{\frac{(n-j)(n-j+2)}{4}-\frac{j(j-1)}{2}}}{(-q;q)_n(q^{n-j+1};q^2)_{\frac{j}{2}}}
, &\text{if $n\equiv j\equiv 0\pmod 2,$}\\[15pt]
\displaystyle (-1)^{\frac{n-j}{2}-1}{n-1\brack \frac{n-1}{2}}_{q^2}
\frac{(q^{n+2};q^2)_{\frac{j-1}{2}}\, q^{\frac{(n-j)(n-j+2)}{4}-\frac{j(j-1)}{2}}}{(-q;q)_{n-1}(q^{n-j+1};q^2)_{\frac{j-1}{2}}}
, &\text{if $n\equiv j\equiv 1\pmod 2,$}\\[15pt]
0, &\text{otherwise.}  \label{eq:sum}
\end{cases}
\end{align}
\end{lem}
\pf
Replacing $k$ by $k+j$, we can write the left-hand side of \eqref{eq:sum} as
\begin{align}
&\hskip -2mm
\sum_{k=0}^{n-j} (-1)^{k+j} {n+k+j\brack n}{n-j\brack k}\frac{q^{\frac{(k+j)(k+j+1)}{2}-n(k+j)}}{(-q;q)_{k+j}} \notag\\
&=\sum_{k=0}^{n-j}(-1)^{k+j}\frac{(q;q)_{n+k+j} (q;q)_{n-j} q^{\frac{(k+j)(k+j+1)}{2}-n(k+j)}}
{(q;q)_{n}(q;q)_{k+j}(q;q)_{k} (q;q)_{n-j-k}(-q;q)_{k+j}} \notag\\
&=\sum_{k=0}^{n-j} (-1)^{j} \frac{(q;q)_{n+j} (q^{n+j+1};q)_{k} (q^{-n+j};q)_{k} q^{\frac{j(j+1)}{2}-nj}}
{(q;q)_{n}(q;q)_{j} (q^{j+1};q)_{k} (q;q)_{k}(-q;q)_{j} (-q^{j+1};q)_{k}} \notag \\
&=(-1)^j{n+j\brack j}\frac{q^{\frac{j(j+1)}{2}-nj}}{(-q;q)_j}
\sum_{k=0}^{n-j} \frac{(q^{n+j+1};q)_{k} (q^{-n+j};q)_{k} q^{k}}{(q;q)_{k} (q^{2j+2};q^2)_{k}}, \label{eq:ktok+j}
\end{align}
where we have used the relation
$$
\frac{(q;q)_{n-j}}{(q;q)_{n-j-k}}=(-1)^k (q^{-n+j};q)_k q^{\frac{k(2n-2j-k+1)}{2}}.
$$

Taking $a=q^{n+j+1}$ and $b=q^{-n+j}$ in  Andrews' $q$-analogue of Gauss' $_2F_1(-1)$ sum
 (see \cite{Andrews73,Andrews74} or \cite[Appendix (II.11)]{GR}):
\begin{align}\label{eq:andrews}
\sum_{k=0}^\infty \frac{(a;q)_k(b;q)_k q^{\frac{k(k+1)}{2}}}{(q;q)_k(abq;q^2)_k}
=\frac{(aq;q^2)_\infty (bq;q^2)_\infty}{(q;q^2)_\infty (abq;q^2)_\infty},
\end{align}
where $(x;q)_\infty=\lim_{n\to\infty}(x;q)_n$,
we have
\begin{align}
\sum_{k=0}^{n-j} \frac{(q^{n+j+1};q)_{k} (q^{-n+j};q)_{k}q^{\frac{k(k+1)}{2}}}
{(q;q)_{k} (q^{2j+2};q^2)_{k}}
&=\frac{(q^{n+j+2};q^2)_\infty (q^{-n+j+1};q^2)_\infty}{(q;q^2)_\infty (q^{2j+2};q^2)_\infty} \notag \\
&=\begin{cases}
\displaystyle \frac{(q^{-n+j+1};q^2)_{\frac{n-j}{2}}}{(q^{2j+2};q^2)_{\frac{n-j}{2}}}, &\text{if $n\equiv j\pmod 2,$}\\[10pt]
0, &\text{otherwise.}  \label{eq:sum-terminate}
\end{cases}
\end{align}
Replacing $q$ by $q^{-1}$ in \eqref{eq:sum-terminate} and noticing that
$(q^{-m};q^{-1})_k=(-1)^k q^{-mk-\frac{k(k-1)}{2}}(q^m;q)_k$, we get
\begin{align}
\sum_{k=0}^{n-j} \frac{(q^{n+j+1};q)_{k} (q^{-n+j};q)_{k}q^{k}}
{(q;q)_{k} (q^{2j+2};q^2)_{k}}
=\begin{cases}
\displaystyle \frac{(q^{-n+j+1};q^2)_{\frac{n-j}{2}}\,q^{(n+j+1)(n-j)}}{(q^{2j+2};q^2)_{\frac{n-j}{2}}}, &\text{if $n\equiv j\pmod 2,$}\\[10pt]
0, &\text{otherwise.}  \label{eq:sum-terminate-2}
\end{cases}
\end{align}
Substituting  \eqref{eq:sum-terminate-2} into \eqref{eq:ktok+j} and making some simplifications,
we obtain the desired identity \eqref{eq:sum}. \qed

\begin{lem}\label{lem:new-legendre}
Let $n$ be a positive integer and
\begin{align}
F_n(x,q)=\sum_{k=0}^n (-1)^k {n\brack k}{n+k\brack k}\frac{(x;q)_k q^{\frac{k(k+1)}{2}-nk}}{(-q;q)_k}.
\label{eq:new-legendre}
\end{align}
Then
\begin{align}
F_n(x,q)=(-1)^n F_n(-x,q). \label{eq:fnxq-sym}
\end{align}
\end{lem}

\pf By the $q$-binomial theorem \eqref{eq:qbinomial}, the coefficient of $x^j$ ($0\leqslant j\leqslant n$)
in the right-hand side of \eqref{eq:new-legendre} is given by
\begin{align*}
&\hskip -2mm
q^{\frac{j(j-1)}{2}}\sum_{k=j}^{n}(-1)^{k-j} {n\brack k} {n+k\brack k}
{k\brack j}\frac{q^{\frac{k(k+1)}{2}-nk}}{(-q;q)_k} \\
&=q^{\frac{j(j-1)}{2}}{n\brack j}\sum_{k=j}^{n}(-1)^{k-j} {n+k\brack k}{n-j\brack k-j}\frac{q^{\frac{k(k+1)}{2}-nk}}{(-q;q)_k}
\end{align*}
which, by Lemma \ref{lem:andrews}, is equal to $0$ if $n-j\equiv 1\pmod 2$. This proves \eqref{eq:fnxq-sym}.
\qed

\noindent{\it Proof of Theorem {\rm\ref{thm:p-121}.}}
Note that, for $k>\frac{p-1}{2}$, there holds $(q;q^2)_k^2\equiv 0\pmod{[p]^2}$.
By Lemma \ref{lem:one}, we have
\begin{align*}
\sum_{k=0}^{p-1}\frac{(q;q^2)_k^2 (x;q^2)_k q^{2k}}{(q^2;q^2)_k^2 (-q^2;q^2)_k}
&\equiv
\sum_{k=0}^{\frac{p-1}{2}}(-1)^k {\frac{p-1}{2}\brack k}_{q^2} {\frac{p-1}{2}+k\brack k}_{q^2}
\frac{(x;q^2)_k q^{k^2+2k-kp}}{(-q^2;q^2)_k} \pmod{[p]^2}.
\end{align*}
The proof then follows from Lemma \ref{lem:new-legendre}. \qed

\noindent{\it Remark.}  Another application of
Andrews's $q$-analogue of Gauss's $_2F_1(-1)$ sum \eqref{eq:andrews} to
$q$-congruences
can be found in \cite{GZ}.
\section{Proofs of Theorems \ref{thm:2mp2}, \ref{thm:2mp2-new} and \ref{thm:main2} }

\noindent{\it Proof of Theorems {\rm\ref{thm:2mp2}}.} 
Since $q^{r}\equiv q^{r-p}\pmod{[p]}$, we may assume that $1\leqslant r\leqslant p$.
When $m=1$, we have
\begin{align*}
\sum_{k=0}^{p-1}\frac{(q^r;q)_k (q^{1-r};q)_k x^k}{(q;q)_k^2}
&=\sum_{k=0}^{r-1} {r-1\brack k}
{r-1+k\brack k} (-x)^k q^{\frac{k(k-1)}{2}-k(r-1)},
\end{align*}
and
\begin{align*}
\sum_{k=0}^{p-1}\frac{(q^r;q)_k (q^{1-r};q)_k (x;q)_k q^{k}}{(q;q)_k^2}
=\sum_{k=0}^{r-1} {r-1\brack k}
{r-1+k\brack k} (-1)^k (x;q)_k q^{\frac{k(k+1)}{2}-k(r-1)}.
\end{align*}
The proof then follows from \eqref{eq:pnxq-1} and \eqref{eq:pnxq-2}
with $P_{n}(x|q)$ replaced by $P_{r-1}(xq^{-1}|q)$.

When $m\geqslant 2$, let
$$
s=\frac{m\langle -\frac{r}{m}\rangle_p+r}{p}.
$$
Then $s$ is a positive integer, $m|ps-r$, and so
\begin{align}
\frac{(q^r;q^m)_k (q^{m-r};q^m)_k}{(q^m;q^m)_k^2}
&=\prod_{j=1}^k \frac{(1-q^{mj-r})(1-q^{mj+r-m})}{(1-q^{mj})^2} \notag\\
&\equiv (-1)^k\prod_{j=1}^k \frac{(1-q^{ps+mj-r})(1-q^{ps-mj-r+m})q^{mj+r-m}}{(1-q^{mj})^2} \notag  \\
&= (-1)^k {\frac{ps-r}{m}\brack k}_{q^m} {\frac{ps-r}{m}+k\brack k}_{q^m} q^{\frac{mk(k-1)}{2}+kr} \notag\\
&\equiv(-1)^k {\langle -\frac{r}{m}\rangle_p \brack k}_{q^m}
 {\langle -\frac{r}{m}\rangle_p+k\brack k}_{q^m} q^{\frac{mk(k-1)}{2}-k(ps-r)}
\pmod{[p]}.  \label{eq:frac-bino-1}
\end{align}
It follows that
\begin{align*}
&\hskip -2mm \sum_{k=0}^{p-1}\frac{(q^r;q^m)_k (q^{m-r};q^m)_k x^k}{(q^m;q^m)_k^2} \notag\\
&\equiv \sum_{k=0}^{\langle -\frac{r}{m}\rangle_p} {\langle -\frac{r}{m}\rangle_p\brack k}_{q^m}
{\langle -\frac{r}{m}\rangle_p+k\brack k}_{q^m} (-x)^k
q^{\frac{mk(k-1)}{2}-mk\langle -\frac{r}{m}\rangle_p}\pmod{[p]},
\end{align*}
and
\begin{align*}
&\hskip -2mm
\sum_{k=0}^{p-1}\frac{(q^r;q^m)_k (q^{m-r};q^m)_k (x;q^m)_k q^{mk}}{(q^m;q^m)_k^2} \notag\\
&\equiv
\sum_{k=0}^{\langle -\frac{r}{m}\rangle_p} {\langle -\frac{r}{m}\rangle_p\brack k}_{q^m}
{\langle -\frac{r}{m}\rangle_p+k\brack k}_{q^m} (-1)^k (x;q^m)_k q^{\frac{mk(k+1)}{2}-mk\langle -\frac{r}{m}\rangle_p}
\pmod{[p]}.
\end{align*}
The proof of \eqref{eq:2mp2} then follows from the two expressions \eqref{eq:pnxq-1} and \eqref{eq:pnxq-2}
for $P_{\langle -\frac{r}{m}\rangle_p}(xq^{-m}|q^m)$. Moreover,
if $p\equiv \pm 1\pmod{m}$, then $\frac{r(m-r)(1-p^2)}{2m}$ is an integer and
$$
\frac{-m\langle -\frac{r}{m}\rangle_p(\langle -\frac{r}{m}\rangle_p+1)}{2}
\equiv \frac{r(m-r)(1-p^2)}{2m}\pmod{p}.
$$
This proves \eqref{eq:2mp2-2}.
\qed

\noindent{\it Proof of Theorems {\rm\ref{thm:2mp2-new}}.} Apply \eqref{eq:pnxq-3}
to $P_{\langle -\frac{r}{m}\rangle_p}(xq^{-m}|q^m)$.  \qed

\noindent{\it Proof of Theorems {\rm\ref{thm:main2}}.} Similarly as before,
we have
\begin{align*}
&\hskip -3mm P_{p-1,m,r}(q,x) \\
&=\sum_{k=0}^{p-1}\frac{(q^r;q^m)_k (q^{m-r};q^m)_k (x;q^m)_k q^{mk}}{(q^m;q^m)_k^2 (-q^m;q^m)_k} \\
&\equiv \sum_{k=0}^{\langle -\frac{r}{m}\rangle_p}(-1)^k  {\langle -\frac{r}{m}\rangle_p\brack k}_{q^m}
{\langle -\frac{r}{m}\rangle_p+k\brack k}_{q^m}
\frac{(x;q^m)_k q^{\frac{mk(k+1)}{2}-mk\langle -\frac{r}{m}\rangle_p}}{(-q^m;q^m)_k}\pmod{[p]}.
\end{align*}
The proof then follows directly from Lemma \ref{lem:new-legendre}.
\qed

\section{Concluding remarks and open problems}
We have the following two stronger conjectural results for Theorems \ref{thm:2mp2}.
\begin{conj}\label{conj:two}
Let $p$ be an odd prime and $m$, $r$ two positive integers with $p\nmid m$ and $m\nmid r$.
Then there exists a unique integer $f_{p,m,r}$ such that
\begin{align*}
&\hskip -2mm \sum_{k=0}^{p-1}\frac{(q^r;q^m)_k (q^{m-r};q^m)_k x^k}{(q^m;q^m)_k^2} \notag\\
&\equiv (-1)^{\langle -\frac{r}{m}\rangle_p}q^{f_{p,m,r}}
\sum_{k=0}^{p-1}\frac{(q^r;q^m)_k (q^{m-r};q^m)_k (x;q^m)_k q^{mk}}{(q^m;q^m)_k^2} \pmod{[p]^2}.
\end{align*}
Furthermore, the numbers $f_{p,m,r}$ satisfy the following recurrence relation:
\begin{align*}
f_{p,m,m+r}
=\begin{cases}
-f_{p,m,r},&\text{if $r\equiv 0\pmod{p},$}\\[5pt]
f_{p,m,r}-r, &\text{otherwise.}
\end{cases}
\end{align*}
\end{conj}
Here are some values of $f_{p,m,r}$:
{\footnotesize
\begin{align*}
&f_{7,2,1}=-12,\ f_{7,2,3}=-13,\ f_{7,2,5}=-16,\  f_{7,2,7}=-21,\
f_{7,2,9}=21,\ f_{7,2,11}=12,\ f_{7,2,13}=1,\   \\
&f_{7,2,15}=-12,\ f_{7,2,17}=-27,\  f_{7,2,19}=-44,\  f_{7,2,21}=-63,\
f_{7,2,23}=63,\ f_{7,2,25}=40,\   \\
&f_{3,5,1}=-5,\ f_{3,5,2}=-3,\ f_{3,5,6}=-6,\  f_{3,5,7}=-5,\
f_{3,5,8}=3,\ f_{3,5,9}=-9,  \\
&f_{7,5,1}=-29,\ f_{7,5,2}=-19,\ f_{7,5,6}=-30,\  f_{7,5,7}=-21,\
f_{7,5,8}=-22,\ f_{7,5,9}=-33, \\
&f_{11,7,1}=-86,\ f_{11,7,2}=-103,\ f_{11,7,3}=-51,\ f_{11,7,8}=-87,\
f_{11,7,9}=-105,  f_{11,7,10}=-54.
\end{align*}}  
\begin{conj}\label{conj:one}
Let $p$ be an odd prime and $m,r$ two positive integers with $r<m$ and $p\equiv \pm 1\pmod{m}$. Then
\begin{align*}
f_{p,m,r}=\frac{r(m-r)(1-p^2)}{2m}.
\end{align*}
\end{conj}

Note that Conjecture~\ref{conj:two} is a $q$-analogue of \eqref{eq:sun-lengendre-2} while
Theorem \ref{thm:main2} is a partial $q$-analogue of \eqref{eq:sun-lengendre-1}, of which
we speculate  the following complete $q$-analogue.
\begin{conj}
Let $p$ be an odd prime and $m$, $r$ two positive integers with $p\nmid m$. Then
\begin{align*}
P_{p-1,m,r}(q,x)
\equiv (-1)^{\langle -\frac{r}{m}\rangle_p} P_{p-1,m,r}(q,-x) \pmod{[p]^2}.
\end{align*}
\end{conj}

There are some similar congruences in the literature. For example, van Hamme \cite{Hamme}
proved the following variant of a conjecture of Beukers \cite{Beukers}:
\begin{align}
\sum_{k=0}^{p-1}\frac{{2k\choose k}^3}{64^k} &\equiv 0
\pmod{p^2},\quad\text{for}\quad p\equiv 3\pmod 4. \label{eq:four-1}
\end{align}
Recently, the authors have obtained a nice $q$-analogue of (7.1),
which will appear in a forthcoming paper.

\vskip 5mm
\noindent{\bf Acknowledgments.} The first author was partially
supported by the Fundamental Research Funds for the Central Universities and
the National Natural Science Foundation of China (grant 11371144).

\end{document}